\documentclass[12pt]{article}
\usepackage[cp866]{inputenc}
\usepackage{citehack}
\usepackage{amsmath}
\usepackage{amsfonts}
\usepackage{amssymb}
\usepackage{longtable}
\newtheorem{theorem}{Theorem}[section]
\newtheorem{prop}{Proposition}[section]
\newtheorem{remark}{Remark}[section]

\newtheorem{lem}{Lemma}[section]

\newtheorem{tab}{Table}
\newtheorem{corol}{Corollary}[section]
\newtheorem{hyp}{Hypothesis}[section]

\def\spa{\mathop\text{{\rm span}}\nolimits}
\def\Hom{\mathop\text{\rm Hom}\nolimits}
\def\pr{\mathop\text{\rm pr}\nolimits}
\def\tr{\mathop\text{\rm tr}\nolimits}

\def\id{\mathop\text{\rm id}\nolimits}

\def\Real{\mathbb{R}}
\def\Co{\mathbb{C}}

\def\g{\mathfrak{g}}
\def\h{\mathfrak{h}}
\def\f{\mathfrak{f}}
\def\z{\mathfrak{z}}
\def\e{\mathfrak{e}}

\def\so{\mathfrak{so}}
\def\osp{\mathfrak{osp}}

\def\sl{\mathfrak{sl}}
\def\psl{\mathfrak{psl}}
\def\sp{\mathfrak{sp}}
\def\gl{\mathfrak{gl}}

\def\hol{\mathfrak{hol}}
\def\q{\mathfrak{q}}

\def\spe{\mathfrak{spe}}
\def\vect{\mathfrak{vect}}
\def\svect{\mathfrak{svect}}
\def\spin{\mathfrak{spin}}

\def\p{\partial}

\def\na{\nabla}

\def\ep{\epsilon}

\def\T{{\cal T}}
\def\M{{\cal M}}

\def\R{{\cal R}}
\def\RR{\bar{\cal R}}
\def\P{{\cal P}}

\title{Irreducible complex skew-Berger algebras}

\author{Anton S. Galaev\footnote{Supported from the Basic Research Center no.
LC505 (Eduard \v{C}ech Center for Algebra and Geometry) of Ministry of Education, Youth and Sport of Czech Republic.} }

\begin{document}

\maketitle\vskip-50ex {\renewcommand{\abstractname}{Abstract}\begin{abstract} Irreducible  skew-Berger algebras
$\g\subset\gl(n,\Co)$, i.e. algebras spanned by the images of the linear maps $R:\odot^2\Co^n\to\g$  satisfying the
Bianchi identity, are classified. These Lie algebras can be interpreted as  irreducible complex Berger superalgebras
contained in $\gl(0|n,\Co)$.

\end{abstract}


{\bf Mathematical subject codes:} 58A50, 53C29

\tableofcontents

\section{Introduction}

The classification of irreducible holonomy algebras of linear torsion-free connections is well known
\cite{Bryant1,Bryant2,M-Sch,Sch}. The first step to this classification was to find candidates to these algebras, these
candidates are called Berger algebras. These algebras $\g\subset\gl(n,\Real)$ are spanned by the images of the linear
maps $R:\Lambda^2\Real^n\to\g$ satisfying the Bianchi identity.

Recently in \cite{Gal} holonomy algebras of connections on supermanifolds were introduced. The natural problem is to
classify irreducible holonomy algebras of linear torsion-free connections on supermanifolds. In \cite{Gal} were defined
Berger superalgebras $\g\subset\gl(m|n,\Real)$, which are the  generalization of the usual Berger algebras, since
$\gl(m|0,\Real)=\gl(m,\Real)$ and Berger superalgebras are the same as Berger algebras in this case. In the present paper
we study the mirror case to the classical one: we classify irreducible complex Berger superalgebras contained in
$\gl(0|n,\Co)$. These Lie superalgebras are the same as irreducible skew-Berger algebras $\g\subset\gl(n,\Co)$, i.e.
algebras spanned by the images of the linear maps $R:\odot^2\Co^n\to\g$ (the skew-curvature tensors) satisfying the Bianchi identity. The reduction
to the real skew-Berger algebras is a standard procedure, see Proposition \ref{Real-Complex} below.

The paper has the following structure. In Section \ref{Spril} we give the necessary preliminaries. In Section \ref{Smain}
we formulate the Main Theorem, where we classify irreducible skew-Berger algebras $\g\subset\gl(n,\Co)$. In Section
\ref{symmetric}  we classify irreducible subalgebras $\g\subset\gl(n,\Co)$ admitting skew-curvature tensors $R$  such that
$R(\Co^n,\Co^n)=\g$ and $\g$ annihilats $R$. This classification immediately follows from the classification of simple
complex Lie superalgebras. Using this list, we obtain examples of skew-Berger algebras. In Section \ref{Sskew-prolongs} we
classify irreducible subalgebras $\g\subset\gl(n,\Co)$ with non-trivial first skew-prolongations
$\g^{[1]}=\{\varphi\in\Hom(\Co^n,\g)|\varphi(x)y=-\varphi(y)x \text{ for all } x,y\in\Co^n\}$ and get examples of
skew-Berger algebras. In Section \ref{Sproof} we finish the proof of the Main Theorem. In Section \ref{Outlook} we explain
how this classification can be used to study some classes of Berger superalgebras $\g\subset\osp(n|2m,\Co)$.

The methods of the paper are mostly taken from \cite{Sch}. In fact, a big number of results of \cite{Sch} can be applied
to our case without change. In the same time some particular cases required new ideas.

{\bf Acknowledgements.} I would like to thank L.~Schwachh\"ofer for useful discussions.
 I am grateful to D.~A.~Leites for the communications concerning the skew-prolongations.
I thank the Department of Mathematics and Statistics of the Masaryk University for the excellent atmosphere for work.

\section{Preliminaries}\label{Spril}

First we give several definitions and facts from \cite{Gal}.

Let $V=V_{\bar 0}\oplus V_{\bar 1}$ be a real or complex vector superspace and  $\g\subset\gl(V)$ a supersubalgebra. {\it
The space of algebraic curvature tensors of type} $\g$ is the vector superspace $\R(\g)=\R(\g)_{\bar 0}\oplus\R(\g)_{\bar
1},$ where $$\R(\g)=\left\{R\in \Lambda^2 V^*
\otimes\g\left|\begin{matrix}R(X,Y)Z+(-1)^{|X|(|Y|+|Z|)}R(Y,Z)X+(-1)^{|Z|(|X|+|Y|)}R(Z,X)Y=0\\\text{for all homogeneous }
X,Y,Z\in V
\end{matrix}\right\}\right..$$ Here $|\cdot|\in\mathbb{Z}_2$ denotes the parity. The identity that satisfy the elements $R\in\RR(\g)$ is called {\it the Bianchi superidentity}.
Obviously, $\R(\g)$ is a $\g$-module with respect to the action \begin{equation}\label{R_A}A\cdot R=R_A,\quad
R_A(X,Y)=[A,R(X,Y)]-(-1)^{|A||R|}R(AX,Y)-(-1)^{|A|(|R|+|X|)}R(X,AY),\end{equation} where $A\in\g$, $R\in\R(\g)$ and
$X,Y\in V$ are homogeneous.

If $\M$ is a supermanifold and $\na$ is a linear torsion-free connection on the tangent sheaf $\T_\M$ with the holonomy
algebra $\hol(\na)_x$ at some point $x$, then for the covariant derivatives of the curvature tensor  we have
$(\na^r_{Y_r,...,Y_1}R)_x\in\R(\hol(\na)_x)$ for all $r\geq 0$ and tangent vectors  $Y_1,...,Y_r\in T_xM$. Moreover,
$|(\na^r_{Y_r,...,Y_1}R)_x|=|Y_1|+\cdots+|Y_r|$, whenever $Y_1,...,Y_r$ are homogeneous.

Define the vector supersubspace $$L(\R(\g))=\spa\{R(X,Y)|R\in\R(\g),\,\,X,Y\in V\}\subset \g.$$ From \eqref{R_A} it
follows that $L(\R(\g))$ is an ideal in $\g$. We call a supersubalgebra $\g\subset\gl(V)$ {\it a Berger superalgebra} if
$L(\R(\g))=\g$.

If $V$ is a vector space, which can be considered as a vector superspace with the trivial odd part,  then
$\g\subset\gl(V)$ is a usual Lie algebra, which can be considered as a Lie superalgebra with the trivial odd part. Berger
superalgebras in this case are the same as the usual Berger algebras.

\begin{prop}  \cite{Gal} Let $\M$ be a supermanifold of dimension $n|m$ with a linear torsion-free connection $\na$.
Then its holonomy algebra $\hol(\na)\subset\gl(n|m,\Real)$ is a Berger superalgebra.\end{prop}

Thus real Berger superalgebras are candidates to the holonomy algebras of linear torsion-free connections on
supermanifolds. The classification of irreducible complex and real Berger algebras is well known
\cite{Bryant1,Bryant2,M-Sch,Sch}.

Consider the vector superspace $$\R^\na(\g)= \left\{S\in V^*
\otimes\R(g)\left|\begin{matrix}S_X(Y,Z)+(-1)^{|X|(|Y|+|Z|)}S_Y(Z,X)+(-1)^{|Z|(|X|+|Y|)}S_Z(X,Y)=0\\ \text{for all
homogeneous } X,Y,Z\in V
\end{matrix}\right\}\right..$$

If $\M$ is a supermanifold and $\na$ is a linear torsion-free connection on $\T_\M$, then
$(\na^r_{Y_r,...,Y_2,\cdot}R)_x\in\R^\na(\hol(\na)_x)$ for all $r\geq 1$ and $Y_2,...,Y_r\in T_xM$. Moreover,
$|(\na^r_{Y_r,...,Y_2,\cdot}R)_x|=|Y_2|+\cdots+|Y_r|$, whenever $Y_2,...,Y_r$ are homogeneous.

A Berger superalgebra $\g$ is called {\it symmetric} if $\R^\na(\g)=0$. This is a generalization of the usual symmetric
Berger algebras, see e.g. \cite{Sch}, and the following is a generalization of the well-known fact about smooth manifolds.

\begin{prop}\cite{Gal} Let $\M$ be a  supermanifold with a torsian free connection $\na$. If $\hol(\na)$ is a symmetric Berger
superalgebra, then $(\M,\na)$ is locally symmetric (i.e. $\nabla R=0$). If $(\M,\na)$ is a locally symmetric superspace,
then its curvature tensor at any point is annihilated by the holonomy algebra at this point and its image coincides with
the holonomy algebra.
\end{prop}

The proof of the following proposition is as in \cite{Sch}.

\begin{prop} Let $\g\subset\gl(V)$ be an irreducible Berger superalgebra. If $\g$ annihilates the module $\R(\g)$,
then $\g$ is a symmetric Berger superalgebra.
\end{prop}

In this paper we consider the case when the vector space $V$ is complex and purely odd, i.e. its even part is trivial. In
this case a supersubalgebra $\g\subset\gl(V)$ is just usual Lie algebra. We may consider the  representation
$\g\subset\gl(\Pi(V))$, where $\Pi$ it the parity changing functor and $\Pi(V)$ becomes a usual vector space.

For a vector space $V$ and a  subalgebra $\g\subset\gl(V)$ define the space of skew-curvature tensors (or just tensors for short) of type $\g$:
$$\RR(\g)=\left\{R\in \odot^2 V^* \otimes\g\left|\begin{matrix}R(X,Y)Z+R(Y,Z)X+R(Z,X)Y=0\\\text{for all  } X,Y,Z\in V
\end{matrix}\right\}\right..$$
Obviously, $\RR(\g)=\R(\g \text{ acting on } \Pi(V))$. A subalgebra $\g\subset\gl(V)$ is called {\it skew-Berger} if
$\g=L(\R(\g)),$ where $$L(\RR(\g))=\spa\{R(X,Y)|R\in\RR(\g),\,\,X,Y\in V\}\subset \g.$$ We see that $\g\subset\gl(V)$ is a
skew-Berger algebra if and only if  $\g\subset\gl(\Pi(V))$ is a Berger superalgebra. Let $\RR^\nabla(\g)=\R^\nabla(\g
\text{ acting on } \Pi(V))$.  A skew-Berger algebra $\g\subset\gl(V)$ is called symmetric if $\RR^\nabla(\g)=0$, i.e.
$\g\subset\gl(\Pi(V))$ is symmetric.

We will use the following fact \begin{equation}\label{vyrazhRR} \RR(\g)=\ker (\p:\odot^2 V^*\otimes\g\to\odot^3V^*\otimes
V),\end{equation} where $\p$ is the symmetrisation map. Note that the map $\p$ is $\g$-equivariant.

Let us explain now how to obtain a classification of irreducible real skew-Berger algebras using the results of this
paper. Let $V$ be a real vector space and $\g\subset\gl(V)$ an irreducible subalgebra. Consider the complexifications
$V_\Co=V\otimes_\Real\Co$ and $\g_\Co=\g\otimes_\Real\Co\subset\gl(V_\Co)$. It is easy to see that
$$\RR(\g_\Co)=\RR(\g)\otimes_\Real \Co\quad\text{ and }\quad \RR^\nabla(\g_\Co)=\RR^\nabla(\g)\otimes_\Real\Co.$$ Recall
that the subalgebra $\g\subset\gl(V)$ is called {\it absolutely irreducible} if $\g_\Co\subset\gl(V_\Co)$ is irreducible
and it is called {\it not absolutely irreducible} otherwise. The last situation appears if and only if there exists a
complex structure $J$ on $V$ commuting with the elements of $\g$. Then $V$ can be considered as a complex vector space and
$\g\subset\gl(V)$ can be considered as a complex irreducible subalgebra. Consider also the natural representation
$i:\g_\Co\to\gl(V)$ in the complex vector space $V$. The following proposition is the analog of Proposition 3.1 from
\cite{Sch}.

\begin{prop}\label{Real-Complex} Let $V$ be a real vector space and $\g\subset\gl(V)$ an irreducible subalgebra.
\begin{description}
\item[1.] If the subalgebra $\g\subset\gl(V)$ is absolutely irreducible, then $\g\subset\gl(V)$ is a skew-Berger algebra
if and only if $\g_\Co\subset\gl(V_\Co)$ is a skew-Berger algebra.

\item[2.] If the subalgebra $\g\subset\gl(V)$ is not absolutely irreducible
and if $(i(\g_\Co))^{[1]}=0$, then  $\g\subset\gl(V)$ is a skew-Berger algebra if and only if $J\g=\g$ and
$\g\subset\gl(V)$ is a complex irreducible skew-Berger algebra.
\end{description}\end{prop}

From this and Proposition 3.1 from \cite{Sch} it follows that if the subalgebra $\g\subset\gl(V)$ is absolutely
irreducible and $\g_\Co\subset\gl(V_\Co)$ is both a skew-Berger and Berger algebra, then the subalgebra $\g\subset\gl(V)$
is a skew-Berger algebra if and only if it is a Berger algebra. Similarly, if the subalgebra $\g\subset\gl(V)$ is not
absolutely irreducible, $(i(\g_\Co))^{[1]}=(i(\g_\Co))^{(1)}=0$ ($(i(\g_\Co))^{(1)}$ denotes the first prolongation), and
$\g\subset\gl(V)$ is both a complex irreducible skew-Berger and Berger algebra, then $\g\subset\gl(V)$ is a real
skew-Berger algebra if and only if it is a real Berger algebra. Thus it is left to consider not absolutely irreducible
real subalgebras $\g\subset\gl(V)$ such that the corresponding representation $i:\g_\Co\to\gl(V)$ in the complex vector
space $V$ is one of the entries 1--10 of Table \ref{tabskew-B} below. This will be done in another paper.

\section{The Main Theorem}\label{Smain}

\begin{theorem}\label{mainTh} Let $V$ be a complex vector space.
The irreducible complex skew-Berger subalgebras $\g\subset\gl(V)$ are exhausted by the representations of Table
\ref{tabskew-B}. The representations 1-8 and 15 have non-trivial first skew-prolongations; the representations 7 with
$\z=0$,  11-14 and 19-22 are symplectic; the representations 8 and 15-18 are orthogonal. If $\g$ admits an element
$\R\in\RR(\g)$ annihilated by $\g$ and such that its image coincides with $\g$, then $\g$ is either 7 with $\z=0$, or 11,
or   19, or 20, or 21, or 22, or 8 with $\g=\sl(n,\Co)$. The representations 19-23 are symmetric skew-Berger algebras. Any
irreducible non-symmetric skew-Berger subalgebras $\g\subset\gl(V)$ coincides with one of the subalgebras 1-18. The
absolutely irreducible real forms of the last two representations can not appear as the holonomy algebras of linear
torsion-free connections on purely odd supermanifolds.
\end{theorem}

\begin{tab}\label{tabskew-B}  Irreducible skew-Berger subalgebras $\g\subset\gl(V)$.
$$
\begin{array}{|c|c|c|c|}\hline &\g& V&\text{restriction}\\\hline\hline
1&\z\oplus\sl(n,\Co)&\Co^n&n\geq 3\\\hline 2&\z\oplus\sl(n,\Co)\oplus\sl(m,\Co)&\Co^n\otimes \Co^m,&n,m\geq 2,\,n\neq
m\\\hline 3&\sl(n,\Co)\oplus\sl(n,\Co)&\Co^n\otimes \Co^n,&n\geq 3\\\hline 4&\sl(n,\Co)&\Lambda^2\Co^n&n\geq 6\\\hline
5&\z\oplus\sl(5,\Co)&\Lambda^2\Co^5&\\\hline 6&\sl(n,\Co)&\odot^2\Co^n&n\geq 3\\\hline
7&\z\oplus\sp(2n,\Co)&\Co^{2n}&n\geq 2\\\hline 8&\g&\g&\g \text{ is a simple Lie algebra}\\\hline\hline
9&\z\oplus\spin(10,\Co)&\Delta^+_{10}=\Co^{16}&\\\hline 10&\f_6&\Co^{27}&\\\hline\hline
11&\sl(2,\Co)\oplus\so(n,\Co)&\Co^2\otimes\Co^{n}&n\geq 3\\\hline 12&\spin(12,\Co)&\Delta^+_{12}=\Co^{32}&\\\hline
13&\sl(6,\Co)&\Lambda^3\Co^6=\Co^{20}&\\\hline 14&\sp(6,\Co)&V_{\pi_3}=\Co^{14}&\\\hline\hline 15&\so(n,\Co)&\Co^n&n\geq
3\\ \hline 16&\g_2&\Co^7&\\\hline 17&\spin(7,\Co)&\Co^8&\\\hline 18&\sl(2,\Co)\oplus\sp(2n,\Co)&\Co^2\otimes\Co^{2n}&n\geq
2\\\hline\hline 19&\sl(2,\Co)&\Co^2&\\\hline 20&\so(n,\Co)\oplus\sp(2m,\Co)&\Co^n\otimes\Co^{2m}&n\geq 3,\,m\geq 2\\\hline
21&\g_2\oplus\sl(2,\Co)&\Co^7\otimes\Co^{2}&\\\hline 22&\spin(7,\Co)\oplus\sl(2,\Co)&\Co^8\otimes\Co^{2}&\\\hline\hline
23&\so(n,\Co)\oplus\sl(m,\Co)&\Co^n\otimes\Co^{m}&n,m\geq 3\\\hline
24&\sp(2n,\Co)\oplus\sl(m,\Co)&\Co^{2n}\otimes\Co^{m}&n\geq 2,\,m\geq 3\\\hline
\end{array}$$
\end{tab}

\begin{remark} From the proof of the theorem it follows that if $\g\subset\gl(V)$ satisfies $\RR(\g)\neq 0$, then $\g\subset\gl(V)$ appears in Table \ref{tabskew-B}, or $\g=\Co\oplus\h$, where $\h\subset\gl(V)$ appears in Table \ref{tabskew-B}.\end{remark}

\begin{remark} We do not find the space $\RR^\nabla(\g)$ for the representations 1-18. If for some of these representation this space is trivial, then absolutely irreducible real forms of this representation can not appear as the holonomy algebras of  linear torsion-free connections on  purely odd supermanifolds.\end{remark}

\begin{remark} The list of representations from Proposition \ref{P3.18} below moustly coincides with the list of representations $\g\subset\gl(V)$ of simple Lie algebras $\g$ such that $\dim \g>\dim V$ \cite{A-V-E}.  We see that for the representations 1-18 of Table \ref{tabskew-B} it holds $\dim \g\geq\dim V$. This proves the following statement.

{\it Let $\g\subset\gl(V)$ be an irreducible skew-Berger algebra. If $\dim \g<\dim V$, then $\g\subset\gl(V)$ is symmetric.}

In fact, nearly the same holds for Berger algebras:

{\it Let $\g\subset\gl(V)$ be an  irreducible Berger algebra. If $\dim \g\leq\dim V$, then $\g\subset\gl(V)$ is symmetric.}

This is the analog of the statement of the Berger holonomy theorem:

{\it Let $G\subset SO(n,\Real)$ be the holonomy group of a Riemannian manifold $(M,g)$. If $G$ does not act transitively on the $n-1$-dimensional sphere, then $(M,g)$ is locally symmetric.}

This formulation follows from the list of possible connected holonomy groups of Riemannian manifolds obtained by M.~Berger in \cite{Berger}. In \cite{Simens} J.~Simens gave a direct proof of this statement. And recently C.~Olmos obtained a more simple and geometric proof of this fact \cite{Olmos}.

\end{remark}

\section{Complex odd symmetric superspaces and the associated skew-Berger algebras}\label{symmetric}

In this section we classify irreducible subalgebras $\g\subset\gl(n,\Co)$ admitting elements $R\in\RR(\g)$ such that
$R(\Co^n,\Co^n)=\g$ and $\g$ annihilats $R$. This classification immediately follows from the classification of simple
complex Lie superalgebras. Then we obtain examples of skew-Berger algebras.

Having such $\g\subset\gl(n,\Co)$ and $R\in\RR(\g)$. Define the Lie superalgebra $\f=\g\oplus\Pi(\Co^n)$ with the
superbrackets $[\xi,\eta]=[\xi,\eta],$ $[\xi,\Pi(x)]=\Pi(\xi x)]$ and $[\Pi(x),\Pi(y)]=R(x,y)$, where $\xi,\eta\in\g$ and
$x,y\in \Co^n$. We get an irreducible infinitesimal symmetric superspace $(\f,\g,\Pi(\Co^n))$ \cite{Cortes2,Serganova}.
The Proposition 1.2.7 from \cite {Kac} implies that $\f$ is a simple Lie superalgebra. In \tab{tabsym} we list simple Lie
superalgebras $\f$ with $\f_{\bar 0}$ acting irreducibly on $\f_{\bar 1}$.

\begin{tab}\label{tabsym} Simple Lie superalgebras $\f$ with $\f_{\bar 0}$ acting irreducibly on $\f_{\bar 1}$.
$$
\begin{array}{|cc|c|c|} \hline
\f&&\f_{\bar 0}&\f_{\bar 1}\\\hline \osp(n|2m,\Co),&n\geq 1,\,\, n\neq 2,\,\,m\geq
1&\so(n,\Co)\oplus\sp(2m,\Co)&\Co^n\otimes\Co^{2m}\\\hline \q(n),&n\geq 3& \sl(n)&\sl(n)\\\hline
F(4)&&\spin(7,\Co)\oplus\sl(2,\Co)&\Co^8\otimes\Co^2\\\hline G(3)&&\g_2\oplus\sl(2,\Co)&\Co^7\otimes\Co^2\\\hline
D(\alpha)&&\sl(2,\Co)\oplus\sl(2,\Co)\oplus\sl(2,\Co)&\Co^2\otimes\Co^2\otimes\Co^2\\\hline
\end{array}$$
\end{tab}

The proof of the following Theorem is similar to the proof of Theorem 3.6 from \cite{Sch}.

\begin{theorem} Let $\g\subset \gl(n,\Co)$ be an irreducible subalgebra. Suppose that there exists an irreducible odd
infinitesimal symmetric superspace
$$((\g\oplus\sl(2,\Co))\oplus\Pi(\Co^n\otimes\Co^2),\g\oplus\sl(2,\Co),\Pi(\Co^n\otimes\Co^2)).$$ Then
$\g\subset\so(n,\Co)$ with respect to some scalar product $g$ on $\Co^n$, there exists a $\g$-equivariant map
$\bar\wedge:\Lambda^2\Co^n\to\g$ satisfying $$(x\bar\wedge y)z+(x\bar\wedge z)y=-2g(y,z)x+g(x,y)z+g(x,z)y$$ for all
$x,y,z\in\Co^n$. Moreover, for any $A\in\g$, the map $R_A:\odot^2\Co^n\to\g$ defined by $$R_A(x,y)=2g(x,y)A+Ax\bar\wedge
y+Ay\bar\wedge x$$ belongs to $\RR(\g)$ and the map $\g\to\RR(\g)$, $A\mapsto R_A$ is injective. In particular, $\g$ is a
skew-Berger algebra.
\end{theorem}

\begin{corol} The following subalgebras are skew-Berger algebras:
$\so(n,\Co)\subset\gl(n,\Co)$ $(n\geq 3)$,  $\sp(2m,\Co)\oplus\sl(2,\Co)\subset\gl(\Co^{2m}\otimes\Co^2)$ $(m\geq 1)$,
$\spin(7,\Co)\subset\gl(8,\Co)$ and $\g_2\subset\gl(7,\Co)$.\end{corol}

{\it Proof.} We get the first algebra using $\osp(n|2,\Co)$ and the fact that $\sp(2,\Co)\simeq \sl(2,\Co)$.
 We get the second algebra using $\osp(4|2m,\Co)$ and the fact that $\so(4,\Co)\simeq\sl(2,\Co)\oplus\sl(2,\Co)$. $\Box$

\begin{remark} Note that $\spin(7,\Co)$ and $\g_2$ are very important  Berger algebras, which turn out to be also
skew-Berger algebras. \end{remark}

\section{Representations with non-trivial first skew-prolongations}\label{Sskew-prolongs}

Let $\g\subset\gl(n,\Co)$ be a subalgebra. Put 
$\g^{[0]}=\g$. Define the $k$-th ($k\geq 1$) skew-prolongation of $\g$ by the rule
$$\g^{[k]}=\{\varphi\in\Hom(\Co^n,\g^{[k-1]})|\varphi(x)y=-\varphi(y)x \text{ for all } x,y\in\Co^n\}.$$


 Let $\g_{-1}$ denote a complex vector superspace and  let  $\g_0\subset\gl(\g_{-1})$ be a supersubalgebra.  {\it The
$k$-th prolongation} ($k\geq 1$) $\g_k$ of $\g_0$ is defined as for the representations of the usual Lie algebras up to
additional signs: $$\g_k=\{\varphi\in\Hom(\g_{-1},\g_{k-1})|\varphi(x)y=(-1)^{|x||y|}\varphi(y)x \text{ for all
homogeneous } x,y\in\g_{-1}\}.$$ Consider {\it the Cartan prolong} $\g_*=\g_*(\g_{-1},\g_0)=\oplus_{k\geq -1}\g_k$. Note
that $\g_*$ has a structure of Lie superalgebra. In \cite{EE,Poletaeva1,Poletaeva2} examples of  irreducible subalgebras
$\g_0\subset\gl(\g_{-1})$ with $\g_1\neq 0$ are given and for the most of them the $(2,2)$-th Spencer cohomology  groups
$H^{2,2}_{\g_0}$ are computed.

It is obvious that for a subalgebra $\g\subset\gl(0|n,\Co)$ its  prolongations coincide with the corresponding
skew-prolongations of the subalgebra $\g\subset\gl(n,\Co)$.

Let $\g\subset\gl(n,\Co)$ be a subalgebra. By analogy with \cite{Sch} we get the following exact sequence
\begin{equation}\label{posled1}0\longrightarrow\g^{[2]}\longrightarrow(\Co^n)^*\otimes \g^{[1]}\longrightarrow\RR(\g)\longrightarrow
H^{2,2}_{\g}\longrightarrow 0,\end{equation} where $H^{2,2}_{\g}$ is the $(2,2)$-th Spencer cohomology group for the
representation $\g\subset\gl(\Pi(\Co^n))$.  The second map in the sequence is given by
\begin{equation}\label{RVposled1}R_{\phi\otimes\alpha}(x,y)=\phi(x)\alpha(y)+\phi(y)\alpha(x).\end{equation}

\begin{theorem}\label{thskewprol}
Let $V$ be a complex vector space. All irreducible subalgebras $\g\subset\gl(V)$ with non-trivial first prolongations and
$\g^{[1]}$, $\g^{[2]}$, $H_\g^{2,2}$ for these subalgebras are listed in Table \ref{tabskewprol}. \end{theorem}

The spaces $(\Co^n\otimes\Lambda^2(\Co^n)^*)_0$ and $(\Co^n\otimes\Lambda^3(\Co^n)^*)_0$ consist of tensors such that the
contraction of the upper index with any down index gives zero. The spaces $H_{\sp(2n,\Co)}^{2,2}$ and
$H_{\sp(2n,\Co)\oplus\Co}^{2,2}$ are given in \cite{Poletaeva2}.

\begin{tab}\label{tabskewprol}  Irreducible subalgebras $\g\subset\gl(n,\Co)$ with $\g^{[1]}\neq 0$.
$$
\begin{array}{|c|c|cc|c|c|c|}\hline &\g& V&&\g^{[1]}&\g^{[2]}&H_\g^{2,2}\\\hline
1&\sl(n,\Co)&\Co^n,&n\geq 3&(\Co^n\otimes\Lambda^2(\Co^n)^*)_0&(\Co^n\otimes\Lambda^3(\Co^n)^*)_0&0 \\\hline
2&\gl(n,\Co)&\Co^n,&n\geq 2&\Co^n\otimes\Lambda^2(\Co^n)^*&\Co^n\otimes\Lambda^3(\Co^n)^*&0 \\\hline
3&\sl(n,\Co)&\odot^2\Co^n,&n\geq 3&\Lambda^2(\Co^n)^*&0&0\\\hline 4&\gl(n,\Co)&\odot^2\Co^n,&n\geq
3&\Lambda^2(\Co^n)^*&0&0\\\hline 5&\sl(n,\Co)&\Lambda^2\Co^n,&n\geq 5&\odot^2(\Co^n)^*&0&0\\\hline
6&\gl(n,\Co)&\Lambda^2\Co^n,&n\geq 5&\odot^2(\Co^n)^*&0&0\text{ if } n\geq 6\\ &&&&&&\Co^5\text{ if } n=5\\\hline
7&\sl(n,\Co)\oplus\sl(m,\Co)\oplus\Co&\Co^n\otimes\Co^m,& n,m\geq 2&V^*&0&0\\ &&&n\neq m&&&\\\hline
8&\sl(n,\Co)\oplus\sl(n,\Co)&\Co^n\otimes\Co^n,& n\geq 3&V^*&0&0\\\hline
9&\sl(n,\Co)\oplus\sl(n,\Co)\oplus\Co&\Co^n\otimes\Co^n,&n\geq 3&V^*&0&0\\\hline 10&\so(n,\Co)&\Co^n,&n\geq 4&\Lambda^3
V^*&\Lambda^4 V^*&0\\\hline 11&\so(n,\Co)\oplus\Co&\Co^n,&n\geq 4&\Lambda^3 V^*&\Lambda^4 V^*&0\\\hline
12&\sp(2n,\Co)\oplus\Co&\Co^{2n}&n\geq 2&V^*&0&\\\hline 13&\g \text{ is simple }&\g& &\Co\id&0&H=?\\\hline
14&\g\oplus\Co,\,\,\g \text{ is simple }&\g& &\Co\id&0&H\\\hline\end{array} $$
\end{tab}

{\bf Proof of Theorem \ref{thskewprol}.} Suppose that for an irreducible subalgebra  $\g\subset\gl(V)$ we have
$\g^{[1]}\neq 0$. Put $\g_{-1}=\Pi(V)$ and $\g_0=\g\subset\gl(\g_{-1})$. We get that $\g_1\neq 0$. It is obvious that the
Cartan prolong $\g_*$ is {\it an irreducible transitive Lie superalgebra with the consistent $\mathbb{Z}$-grading and
$\g_{1}\neq 0$}. This means that $\g_0$ acts irreducibly on $\g_{-1}$, the equality $[a,\g_{-1}]=0$, where $a\in\g_k$,
$k\geq 1$, implies $a=0$, and $(\g_*)_{\bar 0}=\oplus_{k=0}^\infty \g_{2k}$, $(\g_*)_{\bar 1}=\oplus_{k=0}^\infty
\g_{2k-1}$. From \cite[Th. 4]{Kac} it follows that $\g_*$ must coincide with one of the following Lie superalgebras:
\begin{description}\item[I] $\sl(n|m,\Co)$ $(n\neq m,\,\,m,n\geq 2)$, $\psl(n|n,\Co)$ $(n\geq 2)$, $\osp(2|2n,\Co)$
$(n\geq 1)$, $\spe(n,\Co)$ $(n\geq 3)$ with the canonical $\mathbb{Z}$-gradings;
\item[II] $\vect(0|n,\Co)$ $(n\geq 2)$, $\svect(0|n,\Co)$ $(n\geq 3)$, $\h(0|n,\Co)$ $(n\geq 4)$, $\tilde\h(0|n,\Co)$ $(n\geq 4)$ with the canonical $\mathbb{Z}$-gradings;
\item[III] $\tilde\g=\g_{-1}\oplus\g_0\oplus\g_1,$ where $\g_0=\g$ is a simple Lie algebra,
$\g_{-1}=\Pi(\g)$, and $\g_1=\Co$; the non-zero  Lie superbrackets are the following: $[x,y]=[x,y]$,
$[x,\Pi(y)]=\Pi([x,y])$, $[\xi,\Pi(x)]=x$, where $x,y\in\g$ and $\xi\in\Co$;
\item[IV] $\hat\f=\sum_{k=-1}^{\infty}\hat\f_k$ with $\hat\f_0=\f_0\oplus\Co$, $\hat\f_k=\f_k$ for $k\neq 0$ and elements of $\Co$ acting by
the multiplication on $\hat\f_k$ for $k\neq 0$, where $\f$ is of type I, II or III and the center of $\f_0$ is trivial.
\end{description}
The standard $\mathbb{Z}$-gradings of the Lie superalgebras of type I and II are described in \cite{Kac}. The Lie
superalgebras of type I give us the entries 3, 5, 7, 8, 12 of the table; the Lie superalgebras of type II give us 1, 2,
10; the Lie superalgebras of type III give us 13, and the Lie superalgebras of type IV give us all the other entries. We
write $n\geq 5$ for the entries 5 and 6, as $\Lambda^2\Co^3\simeq\Co^2$ and we should consider entries 1 and 2 for $n=3$;
similarly, $\sl(4,\Co)\simeq\so(6,\Co)$ and $\Lambda^2\Co^4\simeq\Co^6$. By an analog reason we assume $n\geq 3$ for the
entries 8 and 9, and $n\geq 2$ for the entry 12. The Spencer cohomology groups $H^{2,2}$ for the entries 1, 2, 3, 5, 7, 8,
10, 12 are computed in \cite{EE,Poletaeva1,Poletaeva2}. The other cohomology groups (except for the entry 13)
 will be computed in the following proposition.

\begin{prop} All the representations of Table \ref{tabskewprol} exsept for the entries 4,
the entries 5 for $n\geq 6$, the entries 9, 11, 12 and 14 are skew-Berger algebras; the representation of
$\sl(n,\Co)\oplus\sl(m,\Co)$ on $\Co^n\otimes\Co^m$ $(n,m\geq 2)$ and $\sp(2n,\Co)\subset\sl(2n,\Co)$ $(n\geq 2)$ are
skew-Berger algebras. \end{prop}

{\bf Proof.} The proof of the fact that the entries 1, 2, 3, 5, 7, 8, 10 of the table and also the representation of
$\sl(n,\Co)\oplus\sl(m,\Co)$ on $\Co^n\otimes\Co^m$ $(n,m\geq 2)$ and $\sp(2n,\Co)\subset\sl(2n,\Co)$ $(n\geq 2)$ are
skew-Berger algebras is similar to the proof of Proposition 12.1 from \cite{Gal}, it follows mostly from the exact
sequence \eqref{posled1}.

\begin{lem} Consider the representation of the Lie algebra $\sl(n,\Co)\oplus\sl(n,\Co)\oplus \Co$ on
$V=\Co^n\otimes\Co^n$ $(n\geq 3)$. Then $\RR(\sl(n,\Co)\oplus\sl(n,\Co)\oplus\Co)=\RR(\sl(n,\Co)\oplus\sl(n,\Co))$ and the
representation of the Lie algebra $\sl(n,\Co)\oplus\sl(n,\Co)\oplus\Co$ on $\Co^n\otimes\Co^n$ $(n\geq 3)$ is not a
skew-Berger algebra.\end{lem}

{\it Proof.} Let $\pi_1,...,\pi_{n-1},\tilde\pi_1,...,\tilde\pi_{n-1}$ denote the fundamental weights of the Lie algebra
$\sl(n,\Co)\oplus\sl(n,\Co)$. Recall that the space $\RR(\sl(n,\Co)\oplus\sl(n,\Co)\oplus\Co)$ can be defined by
\eqref{vyrazhRR}. It can be checked that the $\sl(n,\Co)\oplus\sl(n,\Co)$-module $\odot^2 V^*\simeq\odot^2
V^*\otimes\Co\subset \odot^2 V^*\otimes (\sl(n,\Co)\oplus\sl(n,\Co)\oplus\Co)$ can be decomposed as the direct sum
$V_{\pi_{n-2}+\tilde\pi_{n-2}}\oplus V_{2\pi_{n-1}+2\tilde\pi_{n-1}}$ and each of the decompositions into irreducible
components of the $\sl(n,\Co)\oplus\sl(n,\Co)$-modules  $\odot^2 V^*\otimes(\sl(n,\Co)\oplus\sl(n,\Co))$ and $\odot^3
V^*\otimes V$  contain two copies of the modules $V_{\pi_{n-2}+\tilde\pi_{n-2}}$ and $V_{2\pi_{n-1}+2\tilde\pi_{n-1}}$.
Suppose that $\RR(\sl(n,\Co)\oplus\sl(n,\Co)\oplus\Co)\neq\RR(\sl(n,\Co)\oplus\sl(n,\Co))$. Using the exact sequence
\eqref{posled1}, we get $\RR(\sl(n,\Co)\oplus\sl(n,\Co))\simeq V^*\otimes V^*$ and this
$\sl(n,\Co)\oplus\sl(n,\Co)$-module contains one of each irreducible components $V_{\pi_{n-2}+\tilde\pi_{n-2}}$ and
$V_{2\pi_{n-1}+2\tilde\pi_{n-1}}$. From \eqref{vyrazhRR} it is clear that $\RR(\sl(n,\Co)\oplus\sl(n,\Co)\oplus\Co)$
contains at least one additional component ($V_{\pi_{n-2}+\tilde\pi_{n-2}}$ or  $V_{2\pi_{n-1}+2\tilde\pi_{n-1}}$).
Suppose that $2V_{\pi_{n-2}+\tilde\pi_{n-2}}\subset\RR(\sl(n,\Co)\oplus\sl(n,\Co)\oplus\Co)$. Let $\g_1$ and $\g_2$ denote
the first and the second summands in $\sl(n,\Co)\oplus\sl(n,\Co)$. From the symmetry it follows that each of the
$\sl(n,\Co)\oplus\sl(n,\Co)$-modules $\odot^2 V^*\otimes\g_1$ and $\odot^2 V^*\otimes\g_2$ contains one irreducible
component $V_{\pi_{n-2}+\tilde\pi_{n-2}}$. Denote these modules by $U_1$ and $U_2$. It is clear that any non-zero
$R\in\RR(\sl(n,\Co)\oplus\sl(n,\Co))$ can not take values in one of the Lie algebras $\g_1$ and $\g_2$. Hence $\p|_{U_1}$
and $\p|_{U_2}$ are injective. Moreover, $\p|_{U_1\oplus U_2}$ has the kernel isomorphic to
$V_{\pi_{n-2}+\tilde\pi_{n-2}}$. Denote by $W_1\oplus W_2$, where $W_1$ and $W_2$ are isomorphic to
$V_{\pi_{n-2}+\tilde\pi_{n-2}}$, the submodule $2V_{\pi_{n-2}+\tilde\pi_{n-2}}\subset\odot^3 V^*\otimes V$. We may assume
that $\p|_{U_1}$ and $\p|_{U_2}$ take $U_1$ and $U_2$ isomorphically onto $W_1$. Denote by $U$ the submodule
$V_{\pi_{n-2}+\tilde\pi_{n-2}}\subset \odot^2 V^*\otimes\Co\subset \odot^2 V^*\otimes
(\sl(n,\Co)\oplus\sl(n,\Co)\oplus\Co)$.
It is obviously that $\p|_{U}$ is injective. Let now $S_1+S_2+S_3\in\RR(\sl(n,\Co)\oplus\sl(n,\Co)\oplus\Co)$ and
$S_1+S_2+S_3\not\in\RR(\sl(n,\Co)\oplus\sl(n,\Co))$, where $S_1+S_2+S_3\in U_1\oplus U_2\oplus U$. We have
$\p(S_1+S_2+S_3)=0$. Since $\p(S_1+S_2)\in W_1$ and $\p|_{U_1}:U_1\to W_1$ is an isomorphism, there exists an $S'_1\in
U_1$ such that $\p(S_1+S_2)=\p(S'_1)$. This implies $\p(S'_1+S)=0$, i.e. $S'_1+S$ is a non-zero curvature tensor of type
$\sl(n,\Co)\oplus\sl(n,\Co)\oplus\Co$ taking values in $\g_1\oplus\Co$, which is impossible and we get a contradiction.
The module $V_{2\pi_{n-1}+2\tilde\pi_{n-1}}$ can be considered in the same way. $\Box$

\begin{lem}\label{glLam2} Consider the representation of the Lie algebra $\gl(n,\Co)$ on $V=\Lambda^2\Co^n$. If $n\geq 6$, then $\RR(\gl(n,\Co))=\RR(\sl(n,\Co))$ and $\gl(n,\Co)$ acting on $\Lambda^2\Co^n$ is not a Berger algebra.\end{lem}

{\it Proof.} We have $\odot^2V^*=V_{2\pi_{n-2}}\oplus V_{\pi_{n-4}}$. Suppose that $R\in\RR(\gl(n,\Co))$ and
$R\not\in\RR(\sl(n,\Co))$. Suppose that $R$ has weight $\pi_{n-4}$. Then $R=S+\phi$, where $S\in
\odot^2V^*\otimes\sl(n,\Co)$ and $\phi\in \odot^2V^*\otimes\Co$ have weight $\pi_{n-4}$. Let $e_1,...,e_n$ be the standard
basis of $\Co^n$ and $e_1^*,...,e_n^*$ its dual basis. Assume that $\phi=(e^*_{n-3}\wedge e^*_{n-2})\odot (e^*_{n-1}\wedge
e^*_n)$. Consider the Bianchi identity $$R(e_{n-3}\wedge e_{n-2},e_{n-1}\wedge e_n)e_i\wedge e_j+R(e_{n-1}\wedge
e_n,e_i\wedge e_j)e_{n-3}\wedge e_{n-2}+R(e_i\wedge e_j,e_{n-3}\wedge e_{n-2})e_{n-1}\wedge e_n=0.$$ Note that
$A=S(e_{n-3}\wedge e_{n-2},e_{n-1}\wedge e_n)$ has weight $0$, i.e. it is an element of the Cartan subalgebra
$\mathfrak{t}\subset\sl(n,\Co)$. If $1\leq i\neq j\leq {n-4}$, then $R(e_{n-1}\wedge e_n,e_i\wedge e_j)\in\sl(n,\Co)$ must
has weight $-\epsilon_{n-3}-\epsilon_{n-2}+\epsilon_i+\epsilon_j$ but the Lie algebra $\sl(n,\Co)$ has no such root, i.e.
$R(e_{n-1}\wedge e_n,e_i\wedge e_j)=0$. Similarly, $R(e_{n-3}\wedge e_{n-2},e_i\wedge e_j)=0$. We get that
$(\epsilon_i+\epsilon_j)A=-1$, i.e. $A_{ii}=-\frac{1}{2}$ for $1\leq i\leq {n-4}$ ($A_{ij}$ denote the elements of the
matrix of $A$). By analogy, taking $i={n-3},j={n-2}$ and $i={n-1},j=n$, we get $A_{{n-3}{n-3}}+A_{{n-2}{n-2}}=-1$ and
$A_{{n-1}{n-1}}+A_{nn}=-1$. Thus $\tr A\neq 0$ and we get a contradiction.

The case when $R$ has weight $2\pi_{n-2}$ is similar (we take $\phi=(e^*_{n-1}\wedge e^*_n)\odot (e^*_{n-1}\wedge e^*_n)$). $\Box$

Consider the representation of the Lie algebra $\gl(5,\Co)$ on $V=\Lambda^2\Co^5$. Using the package Mathematica, we find
that $\dim\RR(\gl(n,\Co))=\dim\RR(\sl(n,\Co))+5$. Hence $\gl(5,\Co)$ acting on $\Lambda^2\Co^5$ is  a Berger algebra.
Moreover, $\odot^2 V^*=V_{\pi_1}\oplus V_{2\pi_3}$, hence $\RR(\gl(n,\Co))=\dim\RR(\sl(n,\Co))\oplus U$, where $U$ is
isomorphic to $V_{\pi_1}=\Co^5$.

\begin{lem} Consider the representation of the Lie algebra $\gl(n,\Co)$ on $V=\odot^2\Co^n$,  $n\geq 3$. Then $\RR(\gl(n,\Co))=\RR(\sl(n,\Co))$ and $\gl(n,\Co)$ acting on $\Lambda^2\Co^n$ is not a Berger algebra.
 \end{lem}

{\it Proof.}  The proof of the statement of the lemma for $n\geq 6$ is similar to the proof of the previous lemma (note
that we have $\odot^2V^*=V_{2\pi_{n-1}}\oplus V_{4\pi_{n}}$). To prove the statement for $n=3,4$ and $5$ we use the
package Mathematica. $\Box$

The fact that the subalgebra $\so(n,\Co)\oplus\Co\subset\gl(n,\Co)$ for $n\geq 2$ is not a skew-Berger algebra will be
proved in Proposition \ref{so+C} below. Let $\g$ be a simple Lie algebra and $\g\subset\gl(\g)$ its adjoint
representation. From \eqref{posled1} it follows that $\RR(\g)$ contains a component isomorphic to $\g$ and hence
$\g\subset\gl(\g)$ is a skew-Berger algebra. We do not compute the spencer cohomology for the adjoint representations. Due
to the Cartan-Killing form we have $\g\subset\so(\g)$ and hence $\RR(\g\oplus\Co)=\RR(\g)$ (see Proposition \ref{so+C}
below) and $\g\oplus\Co\subset\gl(\g)$ is not a skew-Berger algebra. The proposition is proved. $\Box$

Consider the representation of the Lie algebra $\sl(n,\Co)\oplus\sl(m,\Co)\oplus \Co$ on the space $V=\Co^n\otimes\Co^m$,
$n,m\geq 2$, $n\neq m$. We get that
 $\RR(\sl(n,\Co)\oplus\sl(m,\Co)\oplus \Co)\simeq V^*\otimes V^*$. To describe this
isomorphism we use the structure of the Lie superbrackets on $\sl(n|m,\Co)$. For $\tau\in V^*\otimes V^*$, the
corresponding curvature tensor is defined by $R_\tau(x_1\otimes x_2,u_1\otimes u_2)=A(x_1\otimes x_2,u_1\otimes
u_2)+B(x_1\otimes x_2,u_1\otimes u_2),$ where $A(x_1\otimes x_2,u_1\otimes u_2)\in\sl(n,\Co)\oplus \Co,$ $B(x_1\otimes
x_2,u_1\otimes u_2)\in\sl(m,\Co)\oplus \Co,$ and for $v_1\in\Co^n$ and $v_2\in\Co^m$ we have $$A(x_1\otimes x_2,u_1\otimes
u_2)v_1=-\tau (x_1,x_2,v_1,u_2)u_1-\tau (u_1,u_2,v_1,x_2)x_1,$$ $$B(x_1\otimes x_2,u_1\otimes u_2)v_2=\tau
(x_1,x_2,u_1,v_2)u_2+\tau (u_1,u_2,x_1,v_2)x_2.$$ In particular, $\tr(R_\tau(x_1\otimes x_2,u_1\otimes
u_2))=(n-m)(\tau(x_1,x_2,u_1,u_2)+\tau(u_1,u_2,x_1,x_2))$. Thus, $\RR(\sl(n,\Co)\oplus\sl(m,\Co))\simeq \Lambda^2V^*$.
Similarly, $\RR(\sl(n,\Co)\oplus\sl(n,\Co))\simeq V^*\otimes V^*$. As we have seen,
$\RR(\sl(n,\Co)\oplus\sl(n,\Co))=\RR(\sl(n,\Co)\oplus\sl(n,\Co)\oplus\Co)$. We will use this later.

\begin{remark} In \cite{Nagy} skew-prolongations of the subalgebras of $\so(n,\Real)$ were considered. In particular it is proved that the only proper irreducible subalgebras of $\so(n,\Real)$ with non-trivial skew-prolongations are exhausted by the adjoint representations of compact simple Lie algebras.
\end{remark}


\section{Proof of Theorem \ref{mainTh}}\label{Sproof}

The following two propositions are analogs of Propositions 3.2 and Lemma 3.5 from \cite{Sch}, respectively.

\begin{prop} Let $\g\subset \sp(2n,\Co)$ be a proper irreducible subalgebra. Then $\RR(\g\oplus\Co)=\RR(\g)$.
In particular, $\g\oplus\Co$ is not a skew-Berger algebra.\end{prop}

\begin{prop}\label{so+C} Let $\g\subset \so(n,\Co)$ be an irreducible subalgebra. Then $\RR(\g\oplus\Co)=\RR(\g)$.
In particular, $\g\oplus\Co$ is not a skew-Berger algebra.\end{prop}

In the following we use the results from \cite{Sch}. In fact, most of these result do not depend on what do we consider:
Berger algebras or skew-Berger algebras (more precisely, it does not matter if we consider a curvature tensor as a map
from $\Lambda^2 V$ or from $\odot^2 V$). Some cases require additional considerations.

Let $V$ be a complex vector space and $\g\subset\gl(V)$ an irreducible subalgebra. Let $\g_s$ and $\z$ denote, respectively, the semi-simple part and the center of $\g$. Then $\g=\g_s\oplus\z$ and either $\z=\Co$ or $\z=0$.
If $\mathfrak t\subset\g$ is a Cartan subalgebra,
let $\mathfrak t_{0}=\mathfrak t\oplus\z$. Denote the set of roots of $\g_s$ by $\Delta$ and the set of weights of the representation $\g\subset\gl(V)$ by $\Phi$.

Let $\Delta_0=\Delta\cup\{0\}$. For each $\alpha\in\Delta$ fix a non-zero $A_\alpha$ in the weight space  $\g_\alpha$, and
let $$\Phi_\alpha=\{\text{weights of } A_\alpha V\}=(\alpha+\Phi)\cap\Phi.$$

A triple $(\lambda_0,\lambda_1,\alpha)$, where $\lambda_0,\lambda_1\in \Phi$ and $\alpha\in\Delta$, is called {\it a
spanning triple} if $$\Phi_\alpha\subset\{\lambda_0+\beta,\lambda_1+\beta|\beta\in\Delta_0\}.$$ A spanning triple is
called {\it extremal} if $\lambda_0$ and $\lambda_1$ are extremal weights.

\begin{prop}\cite[Prop. 3.10]{Sch} Let $\g\subset\gl(V)$ be an irreducible skew-Berger algebra.
Then for every root $\alpha\in\Delta$ there is a spanning triple $(\lambda_0,\lambda_1,\alpha)$, a weight element
$R\in\RR(\g)$ and vectors $x_0,x_1\in V$ of weights $\lambda_0,\lambda_1$ such that $R(x_0,x_1)=A_\alpha$.

In fact, if $R\in\RR(\g)$ is a weight element and if there are weight vectors $x_0,x_1\in V$ of weights
$\lambda_0,\lambda_1$ such that $R(x_0,x_1)=A_\alpha$, then $(\lambda_0,\lambda_1,\alpha)$ is a spanning triple.
\end{prop}

\begin{theorem} Let $\g\subset\gl(V)$ be an irreducible skew-Berger algebra. Then  there is an extrimal spanning triple $(\lambda_0,\lambda_1,\alpha)$.
\end{theorem}

{\bf Proof.} The proof is the same as the proof of Theorem 3.12 from \cite{Sch}. As the last step we need to show that
$\so(n,\Co)$ acting on $(\odot^2\Co^n)_0=\odot^2\Co^n/\Co g,$ where $g$ is the scalar product on $\Co^n$,  is not a
skew-Berger algebra. Consider the inclusion $\so(n,\Co)\subset\sl(n,\Co)$ and the representation of the Lie algebra
$\sl(n,\Co)$ on the vector space $\odot^2\Co^n$. Then, $\RR(\so(n,\Co))\subset\RR(\sl(n,\Co))$. From Section
\ref{Sskew-prolongs} we know that $\RR(\sl(n,\Co))\simeq \odot^2(\Co^n)^*\otimes\Lambda^2(\Co^n)^*$. To describe this
isomorphism we use the structure of the Lie superbrackets of $\spe(n,\Co)$. For
$\tau\in\odot^2(\Co^n)^*\otimes\Lambda^2(\Co^n)^*$, the corresponding curvature tensor is defined by $$R_\tau(x_1\odot
x_2,y_1\odot
y_2)z=-2\big(\tau(x_1,x_2,y_1,z)y_2+\tau(x_1,x_2,y_2,z)y_1+\tau(y_1,y_2,x_1,z)x_2+\tau(y_1,y_2,x_2,z)x_1\big),$$ where
$x_1,x_2,y_1,y_2,z\in\Co^n$. It is not hard to verify that from the condition $R_\tau(x_1\odot x_2,y_1\odot
y_2)\in\so(n,\Co)$ for all $x_1,x_2,y_1,y_2\in\Co^n$ it follows that $R_\tau=0$. Thus, $\R(\so(n,\Co))=0$. $\Box$

Now we will consider the case when $\g_s$ is simple.

\begin{prop}\cite[Prop. 3.18]{Sch}\label{P3.18}
Let $\g\subset\gl(V)$ be an irreducible subalgebra such that $\g_s$ is simple. Suppose that there exists an extremal spanning triple $(\lambda_0,\lambda_1,\alpha)$. Then either the dominant weight is a root, i.e. $\Phi\subset\Delta_0$, or the representation of $\g_s$ on $V$ is conjugated to one of the following:

\vskip0.3cm

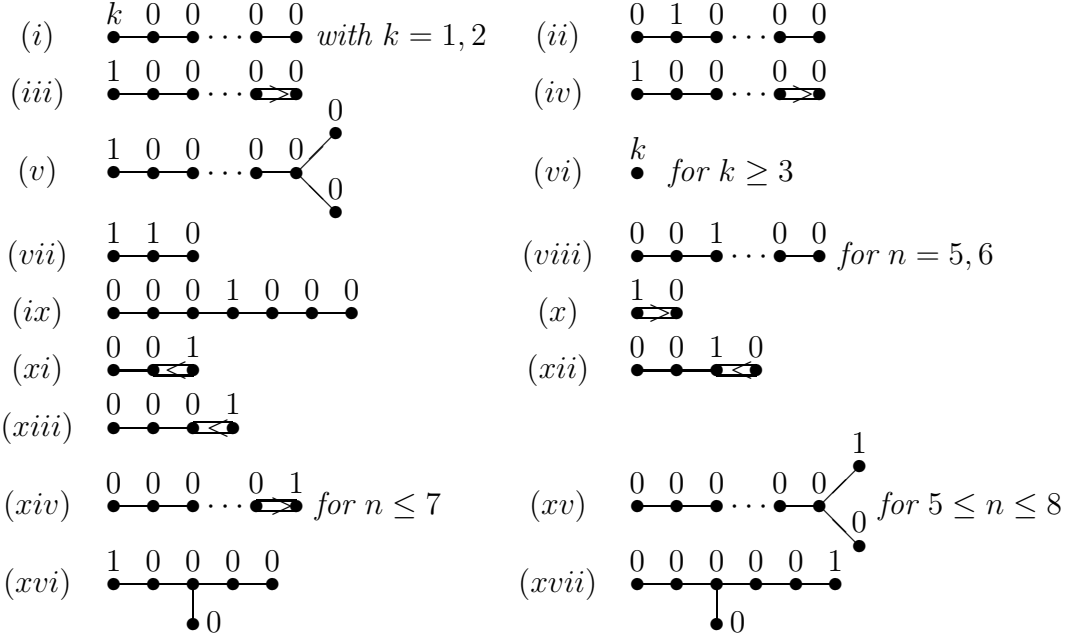
\begin{tabular}{clcl}
$(i)$ & \begin{picture}(75,16)
\put(0,0){$\bullet$}\put(0,8){$k$}
\put(3,3){\line(1,0){15}}\put(15,0){$\bullet$}\put(15,8){$0$}
\put(18,3){\line(1,0){15}}
\put(30,0){$\bullet\cdots\bullet$}\put(30,8){$0$}\put(54,8){$0$}\put(57,3)
{\line(1,0){15}}\put(69,0){$\bullet$}\put(69,8){$0$}\end{picture} with $k=1,2$&
$(ii)$ & \begin{picture}(110,16)\put(0,0){$\bullet$}\put(0,8){$0$}\put(3,3){\line(1,0){15}}\put(15,0){$\bullet$}\put(15,8){$1$}
\put(18,3){\line(1,0){15}}
\put(30,0){$\bullet\cdots\bullet$}\put(30,8){$0$}\put(54,8){$0$}\put(57,3)
{\line(1,0){15}}\put(69,0){$\bullet$}\put(69,8){$0$}\end{picture}\\

$(iii)$ & \begin{picture}(75,16)\put(0,0){$\bullet$}\put(0,8){$1$}\put(3,3){\line(1,0){15}}\put(15,0){$\bullet$}\put(15,8){$0$}
\put(18,3){\line(1,0){15}}
\put(30,0){$\bullet\cdots\bullet$}\put(30,8){$0$}\put(54,8){$0$}\put(57,1)
{\line(1,0){15}}\put(57,5)
{\line(1,0){15}}\put(61,0){$>$}\put(69,0){$\bullet$}\put(69,8){$0$}\end{picture}&
$(iv)$ & \begin{picture}(75,16)\put(0,0){$\bullet$}\put(0,8){$1$}\put(3,3){\line(1,0){15}}\put(15,0){$\bullet$}\put(15,8){$0$}
\put(18,3){\line(1,0){15}}
\put(30,0){$\bullet\cdots\bullet$}\put(30,8){$0$}\put(54,8){$0$}\put(57,1)
{\line(1,0){15}}\put(57,5)
{\line(1,0){15}}\put(61,0){$>$}\put(69,0){$\bullet$}\put(69,8){$0$}\end{picture}\\

$(v)$ & \begin{picture}(90,24)\put(0,0){$\bullet$}\put(0,8){$1$}\put(3,3){\line(1,0){15}}\put(15,0){$\bullet$}\put(15,8){$0$}
\put(18,3){\line(1,0){15}}
\put(30,0){$\bullet\cdots\bullet$}\put(30,8){$0$}\put(54,8){$0$}\put(57,3)
{\line(1,0){15}}\put(69,0){$\bullet$}\put(69,8){$0$}\put(72,3){\line(1,1){15}}\put(84,15){$\bullet$}\put(84,23){$0$}
\put(72,3){\line(1,-1){15}}\put(84,-15){$\bullet$}\put(84,-7){$0$}\end{picture}&
$(vi)$ & \begin{picture}(15,16)\put(0,0){$\bullet$}\put(0,8){$k$} \end{picture}for $k\geq 3$\\

$(vii)$ & \begin{picture}(75,26)\put(0,0){$\bullet$}\put(0,8){$1$}\put(3,3){\line(1,0){15}}\put(15,0){$\bullet$}\put(15,8){$1$}
\put(18,3){\line(1,0){15}}
\put(30,0){$\bullet$}\put(30,8){$0$}
\end{picture}&
$(viii)$ & \begin{picture}(75,26)\put(0,0){$\bullet$}\put(0,8){$0$}\put(3,3){\line(1,0){15}}\put(15,0){$\bullet$}\put(15,8){$0$}
\put(18,3){\line(1,0){15}}
\put(30,0){$\bullet\cdots\bullet$}\put(30,8){$1$}\put(54,8){$0$}\put(57,3)
{\line(1,0){15}}\put(69,0){$\bullet$}\put(69,8){$0$}\end{picture} for $n=5,6$\\

$(ix)$&
\begin{picture}(85,16)\put(0,0){$\bullet$}\put(0,8){$0$}\put(3,3){\line(1,0){15}}\put(15,0){$\bullet$}\put(15,8){$0$}
\put(18,3){\line(1,0){15}}
\put(30,0){$\bullet$}\put(30,8){$0$}
\put(33,3){\line(1,0){15}}\put(45,0){$\bullet$}\put(45,8){$1$}
\put(48,3){\line(1,0){15}}\put(60,0){$\bullet$}\put(60,8){$0$}
\put(63,3){\line(1,0){15}}\put(75,0){$\bullet$}\put(75,8){$0$}
\put(78,3){\line(1,0){15}}\put(90,0){$\bullet$}\put(90,8){$0$}
\end{picture}&
$(x)$&
\begin{picture}(40,16)\put(0,0){$\bullet$}\put(0,8){$1$}\put(3,1){\line(1,0){15}}\put(3,5){\line(1,0){15}}\put(7,0){$>$}
\put(15,0){$\bullet$}\put(15,8){$0$}\end{picture}\\

$(xi)$&
\begin{picture}(65,16)
\put(0,0){$\bullet$}\put(0,8){$0$}
\put(3,3){\line(1,0){15}}\put(15,0){$\bullet$}\put(15,8){$0$}
\put(18,1){\line(1,0){15}}\put(18,5){\line(1,0){15}}\put(22,0){$<$}
\put(30,0){$\bullet$}\put(30,8){$1$}
\end{picture}&
$(xii)$&
\begin{picture}(65,16)
\put(0,0){$\bullet$}\put(0,8){$0$}
\put(3,3){\line(1,0){15}}\put(15,0){$\bullet$}\put(15,8){$0$}
\put(18,3){\line(1,0){15}}\put(30,0){$\bullet$}\put(30,8){$1$}
\put(33,1){\line(1,0){15}}\put(33,5){\line(1,0){15}}\put(38,0){$<$}
\put(45,0){$\bullet$}\put(45,8){$0$}
\end{picture}\\
$(xiii)$&
\begin{picture}(65,16)
\put(0,0){$\bullet$}\put(0,8){$0$}
\put(3,3){\line(1,0){15}}\put(15,0){$\bullet$}\put(15,8){$0$}
\put(18,3){\line(1,0){15}}\put(30,0){$\bullet$}\put(30,8){$0$}
\put(33,1){\line(1,0){15}}\put(33,5){\line(1,0){15}}\put(38,0){$<$}
\put(45,0){$\bullet$}\put(45,8){$1$}
\end{picture}&&\\

$(xiv)$ & \begin{picture}(75,16)
\put(0,0){$\bullet$}\put(0,8){$0$}
\put(3,3){\line(1,0){15}}\put(15,0){$\bullet$}\put(15,8){$0$}
\put(18,3){\line(1,0){15}}
\put(30,0){$\bullet\cdots\bullet$}\put(30,8){$0$}\put(54,8){$0$}
\put(57,1){\line(1,0){15}}\put(57,5){\line(1,0){15}}
\put(62,0){$>$}
\put(69,0){$\bullet$}\put(69,8){$1$}\end{picture} for  $n\leq 7$&

$(xv)$ & \begin{picture}(90,24)\put(0,0){$\bullet$}\put(0,8){$0$}\put(3,3){\line(1,0){15}}\put(15,0){$\bullet$}\put(15,8){$0$}
\put(18,3){\line(1,0){15}}
\put(30,0){$\bullet\cdots\bullet$}\put(30,8){$0$}\put(54,8){$0$}\put(57,3)
{\line(1,0){15}}\put(69,0){$\bullet$}\put(69,8){$0$}\put(72,3){\line(1,1){15}}\put(84,15){$\bullet$}\put(84,23){$1$}
\put(72,3){\line(1,-1){15}}\put(84,-15){$\bullet$}\put(84,-7){$0$}\end{picture} for $5\leq n\leq 8$ \\

$(xvi)$ &
\begin{picture}(85,24)\put(0,0){$\bullet$}\put(0,8){$1$}\put(3,3){\line(1,0){15}}\put(15,0){$\bullet$}\put(15,8){$0$}
\put(18,3){\line(1,0){15}}
\put(30,0){$\bullet$}\put(30,8){$0$}
\put(33,3){\line(1,0){15}}\put(45,0){$\bullet$}\put(45,8){$0$}
\put(48,3){\line(1,0){15}}\put(60,0){$\bullet$}\put(60,8){$0$}
\put(33,3){\line(0,-1){15}}
\put(30,-15){$\bullet$}\put(38,-15){$0$}
\end{picture}&
$(xvii)$&
\begin{picture}(85,24)\put(0,0){$\bullet$}\put(0,8){$0$}\put(3,3){\line(1,0){15}}\put(15,0){$\bullet$}\put(15,8){$0$}
\put(18,3){\line(1,0){15}}
\put(30,0){$\bullet$}\put(30,8){$0$}
\put(33,3){\line(1,0){15}}\put(45,0){$\bullet$}\put(45,8){$0$}
\put(48,3){\line(1,0){15}}\put(60,0){$\bullet$}\put(60,8){$0$}
\put(33,3){\line(0,-1){15}}
\put(30,-15){$\bullet$}\put(38,-15){$0$}
\put(63,3){\line(1,0){15}}\put(75,0){$\bullet$}\put(75,8){$1$}
\end{picture}
\end{tabular}

\vskip0.5cm
\end{prop}


Note that if $\Phi\subset\Delta_0$, then $\g_s\subset\gl(V)$ is one of the following:

$(xviii)$ $\sp(2n,\Co)$ acting on $(\Lambda^2\Co^{2n})_0=\Lambda^2\Co^{2n}/\Co\Omega$, where $\Omega$ is the symplectic form on $\Co^{2n}$;

$(xix)$ $\f_4$ acting on $V^{\f_4}_{\pi_1}=\Co^{26}$;

$(xx)$ $\g_2\subset\gl(7,\Co)$;

$(xxi)$ the adjoint representation of a simple Lie algebra.

\vskip0.3cm

We have already discussed the representations $(i)-(v)$, $(xiv)$ for $n=3$, $(xx)$  and $(xxi)$.

Now we consider the remaining representations.

\vskip0.2cm

$(vi)$ ($\gl(2,\Co)\subset\gl(\odot^k\Co^2),\,\,k\leq 3$). The computations using the package Mathematica show that
$\R(\g_s)=1,3 ,0 $ for $k=1,2,3$, respectively, and $\R(\g_s)=6,3,0$ for $k=1,2,3$. Thus the only skew-Berger algebras are
$\g_s$ and $\g$ for $k=1$, and $\g_s$ for $k=2$ (in the last case we get the standard representation of $\so(3,\Co)$).

 \vskip0.2cm

The representations $(vii)$, $(ix)$, $(x)$, $(xii)$, $(xiii)$, $(xv)$, $n=8$, can be dealt exactly in the same way as in
\cite{Sch} (for some of these representations it is proved $\dim\R(\g)\leq 1$; in the same way we get $\dim\RR(\g)\leq 1$,
but these representation does not appear in Table \ref{tabskewprol}, hence $\dim\RR(\g)=0$).

\vskip0.2cm

Recall that the space $\RR(\g)$ can be find from \eqref{vyrazhRR} and the map $\p$ is $\g$-equivariant. Decompose the
$\g$-modules $\odot^2 V^*\otimes\g$ and $\odot^3 V^* \otimes V$ into the  direct sums of irreducible components. If a
component $V_\Lambda$ appears in $\odot^2 V^*\otimes\g$ more times than in $\odot^3 V^* \otimes V$, then the space
$\RR(\g)$ containes a $\g$-submodule isomorphic to $V_\Lambda$. In particular, $\RR(\g)\neq 0$, and if $\g$ is simple,
then it is a skew-Berger algebra. We apply this idea to the following 4 cases.

\vskip0.2cm

$(xv)$ for $n=5$ (the spin representation of the Lie algebra $\g=\so(10,\Co)$ on $V=\Delta_{10}^+)$. Using the package
LiE, we check that $\odot^2 V^*\otimes\g$ containes two copies of $V_{\pi_3}$, while $\odot^3 V^* \otimes V$ containes
only one copy of $V_{\pi_3}$. Note that $\dim V_{\pi_3}=120$. Using the programm Mathematica we find that
$\dim\RR(\spin(10,\Co))=120$. Thus, $\RR(\spin(10,\Co))\simeq V_{\pi_3}$. Since the Lie algebra $\so(10,\Co)$ is simple,
we get that the representation $\spin(10,\Co)\subset\gl(\Delta_{10}^+)$ is a skew-Berger algebra.  Further more, we find
that $\dim\RR(\spin(10,\Co)\oplus\Co)=176$. Consequently, $\spin(10,\Co)\oplus\Co\subset\gl(\Delta_{10}^+)$ is a
skew-Berger algebra.

\vskip0.2cm

$(xv)$ for $n=6$ (the spin representation of the Lie algebra $\g=\so(12,\Co)$ on $V=\Delta_{12}^+)$. Using the package
LiE, we check that $\odot^2 V^*\otimes\g$ containes a submodule isomorphic to $V_{2\pi_1}$, while $\odot^3 V^* \otimes V$
does not contain a submodule isomorphic to $V_{2\pi_1}$. Thus $\spin(12,\Co)\subset\gl(\Delta_{12}^+)$ is a skew-Berger
algebra. Since $\spin(12,\Co)$ is contained in the symplectic Lie algebra,
$\spin(12,\Co)\oplus\Co\subset\gl(\Delta_{12}^+)$ is not a skew-Berger algebra.

\vskip0.2cm

$(xvi)$ ($\e_6\subset \gl(27,\Co)$). Using the package LiE, we check that $\odot^2 V^*\otimes\g$ containes two copies of
$V_{\pi_5}$, while $\odot^3 V^* \otimes V$ containes only one copy of $V_{\pi_5}$. Thus $\e_6$ acting on $\Co^{27}$ is a
skew-Berger algebra. We claim that $\RR(\e_6\oplus\Co)=\RR(\e_6)$, i.e. $\e_6\oplus\Co$ acting on $\Co^{27}$ is not a
skew-Berger algebra. The idea of the proof is similar to the proof of Lemma \ref{glLam2}.  Note that the representation of
$\g$ on $V$ has 27 weights. This means that each weight subspace of $V$ is one-dimensional. For each weight $\lambda$
choose a non-zero weight vector $e_\lambda\in V$. Then the vectors $e_\lambda$ form a basis of $V$. Let $e^*_\lambda$ be
the dual basis. We have $\odot^2V^*=V_{2\pi_6}\oplus V_{\pi_1}$. Suppose that $S+\phi\in\RR(\e_6\oplus\Co)$, $S$ and
$\phi$  have weight $\pi_1$ and $\phi\neq 0$. We may assume that $\phi=e^*_{-\pi_6}\odot e^*_{-\pi_1+\pi_6}$. Note that
$A=S(e_{-\pi_6}, e_{-\pi_1+\pi_6})$ has weight 0, i.e. it is an element of the Cartan subalgebra of $\e_6$. Writing down
the Bianchi identity for the vectors $e_{-\pi_6}$, $e_{-\pi_6}$, $e_{-\pi_1+\pi_6}$, and using the fact that
$\pi_1-2\pi_6$ is not a root of $\e_6$, we get  $R(e_{-\pi_6}, e_{-\pi_1+\pi_6})e_{-\pi_6}=0$. Consequently, $\pi_6(A)=1$.
Similarly, writing the Bianchi identity for the vectors $e_{-\pi_6}$, $e_{-\pi_1+\pi_6}$, $e_{-\pi_1+\pi_6}$ and for the
vectors $e_{-\pi_6}$, $e_{\pi_1}$, $e_{-\pi_1+\pi_6}$, we get that $(-\pi_1+\pi_6)A=-1$ and $\pi_1(A)=-1$, which leads to
a contradiction. The elements $S+\phi\in\RR(\e_6\oplus\Co)$ of weight $2\pi_6$ can be considered in the same way. We
choose $\phi=e^*_{-\pi_6}\odot e^*_{-\pi_6}$ and use the Bianchi identity for the vectors $e_{-\pi_6}$, $e_{-\pi_6}$,
$e_{-\pi_6}$, for the vectors $e_{-\pi_6}$, $e_{-\pi_6}$, $e_{-\pi_3+\pi_4}$, and for the vectors $e_{-\pi_6}$,
$e_{-\pi_6}$, $e_{\pi_3-\pi_4-\pi_6}$.

\vskip0.2cm

$(xi)$ ($\sp(6,\Co)\subset\gl(\Lambda^3\Co^6)$). In this case $\odot^2 V^*\otimes\g$ containes a submodule isomorphic
to $V_{\pi_3}$, while $\odot^3 V^* \otimes V$ does not contain a submodule isomorphic to $V_{\pi_3}$. Hence
$\sp(6,\Co)\subset \gl(14,\Co)$ is a skew-Berger algebra. Since this representation is symplectic, $\g\oplus\Co$ is not a
skew-Berger algebra.

\vskip0.2cm

$(xv)$ for $n=7$ (the spin representation of the Lie algebra $\g_s=\so(14,\Co)$ on $V=\Delta_{14}^+$). The only spanning
triples, up to the action of the Weyl group, are $(\pi_7,\pi_3-\pi_7,\pi_2)$ and $(\pi_7,\pi_1-\pi_7,\pi_2)$ \cite{Sch}.
Consequently, elements of $\RR(\g)$ may have, up to the action of the Weyl group, weights $\pi_2-\pi_3$ and $\pi_2-\pi_1$.
Suppose that $R\in\RR(\g)$ has weight $\pi_2-\pi_3$. Let $x,\,\,y,\,\,z\in V$ be vectors of weights $\pi_7$, $\pi_3-\pi_7$
and $-\pi_2+\pi_3-\pi_7$. Considering the Bianchi identity for these vectors and using the facts that $\pi_3$ and
$\pi_3-2\pi_7$ are not roots, we get $R(y,z)=R(x,z)=0$ and $R(x,y)z=0$. If we suppose that $R(x,y)\neq 0$, then since
$R(x,y)$ is a root vector of weight $\pi_2$, $z$ has weight $-\pi_2+\pi_3-\pi_7$, and $\pi_3-\pi_7$ belongs to the weights
of the representation, we get $R(x,y)z\neq 0$. Thus, $R(x,y)=0$. Similarly, if $x$ and $y$ have weights  $\pi_7$ and
$\pi_1-\pi_7$, respectively, then in the same way we show that $R(x,y)=0$ (choose $z$ of weight $-\pi_2+\pi_5-\pi_7$).
Consider a spanning triple $(\lambda_0,\lambda_1,\pi_2)$. Let $x$ and $y$ have weights $\lambda_0$ and $\lambda_1$,
respectively. Then there is an element $w$ in the Weyl group taking $(\lambda_0,\lambda_1,\pi_2)$ either to
$(\pi_7,\pi_3-\pi_7,\pi_2)$ or to $(\pi_7,\pi_1-\pi_7,\pi_2)$, let us assume the first. Applying $w$ to the system of
positive roots, we get another system of positive roots. In this new system $(\lambda_0,\lambda_1,\pi_2)$ has the same
expression as
 $(\pi_7,\pi_3-\pi_7,\pi_2)$ in the old one. Consequently,
$R(x,y)\in \g_{\pi_2}$ implies $R(x,y)=0$. And there are no weight elements $R$ and $x,y$ such that $0\neq R(x,y)\in
\g_{\pi_2}$. Since $\g_s$ is simple, we conclude that $\RR(\g)=0$. Thus $\spin(14,\Co)\subset\gl(\Delta_{14}^+)$ and
$\spin(14,\Co)\oplus\Co\subset\gl(\Delta_{14}^+)$ are not skew-Berger algebras.

\vskip0.2cm

$(xiv)$ (the spin representation of the Lie algebra $\g_s=\so(2n+1,\Co)$ on $V=\Delta_{2n+1}$, $n\leq 7$). We have
$\spin(2n+1,\Co)\subset\spin(2n+2,\Co)$. From the above we get $\RR(\spin(2n+1,\Co))=0$ for $n=6,7$. We use the package
Mathematica to show that $\RR(\spin(2n+1,\Co))=0$ for $n=4,5$. The case $n=3$ is already considered in Section
\ref{symmetric}; alternatively, using the package Mathematica we find that $\dim\RR(\spin(7,\Co))=126$. If $n=2$, then we
get the standard representation of $\sp(4,\Co)$; for $n=1$, we get the standard representation of $\sl(2,\Co)$.

\vskip0.2cm

$(viii)$ If $n=6$, then using the package Mathematica we show that $\dim\RR(\g_s)=35$; in this case the representation is
symplectic and $\RR(\g)=\RR(\g_s)$. Suppose that n=7. The only spanning triples, up to the action of the Weyl group, are
$(\ep_1+\ep_2+\ep_3,\ep_4+\ep_5+\ep_6,\ep_1-\ep_7)$ and $(\ep_1+\ep_2+\ep_3,\ep_1+\ep_4+\ep_5,\ep_1-\ep_7)$ \cite{Sch}. We
will denote a vector $e_i\wedge e_j\wedge e_k$ by $e_{ijk}$. Suppose that $R\in\RR(\g)$ is a weight vector and
$R(e_{123},e_{456})\in\g_{\ep_1-\ep_7}$, i.e. $R$ has weight $\ep_1$. There is some $a\in \Co$ such that
$R(e_{123},e_{456})=a E_{17}$, where $E_{17}$ is the matrix with $1$ on the position $(1,7)$ and $0$ on the other
positions. Let $A=R(e_{237},e_{456})\in\g$. Then $A$ is an element of the Cartan subalgebra of $\g$. Applying the Bianchi
identity to the vectors $e_{123}$, $e_{456}$ and $e_{237}$, and using that $R(e_{123},e_{237})=0$ (as this element has
weight which is not a root of $\g$), we get $A_{11}+A_{22}+A_{33}+a=0$. Applying the Bianchi identity to $e_{456}$,
$e_{237}$ and one of the vectors $e_{12i}$ ($4\leq i\leq 6$), $e_{237}$,  $e_{147}$, we get  $A_{11}+A_{22}+A_{ii}=0$,
$A_{22}+A_{33}+A_{77}=0$, $A_{11}+A_{44}+A_{77}=0$. Using these conditions and the traceless of $A$, we get that $a=0$.
The spanning triple $(\ep_1+\ep_2+\ep_3,\ep_1+\ep_4+\ep_5,\ep_1-\ep_7)$ can be considered in the same way. As  in the case
$(xv)$ for $n=7$ we conclude that if $R\in\RR(\g)$, $x$ and $y$ are weight elements such that $R(x,y)\in\g_{\ep_1-\ep_7}$,
then $R(x,y)=0$, and $\RR(\g)=0$.

\vskip0.2cm

$(xvii)$:  $\e_7\subset\gl(56,\Co)$ is not a skew-Berger algebra. In \cite{Adams} there is the following description of
this representation. The Lie algebra $\e_7$ admits the following structure of $\mathbb{Z}_2$-graded Lie algebra:
$\e_7=\sl(8,\Co)\oplus\Lambda^4\Co^8$. The representation space $\Co^{56}$ is decomposed into the direct sum
$\Co^{56}=\Lambda^2\Co^8\oplus\Lambda^2(\Co^8)^*$. The  elements of $\sl(8,\Co)$ preserve this decomposition and act
naturally on each component, and the elements of $\Lambda^4\Co^8$ interchange these components. Since
$\e_7\subset\sp(56,\Co)$, any element $R\in\RR(\e_7)$ is a symmetric map $R:\odot^2\Co^{56}\to \e_7\subset
\odot^2\Co^{56}$ and it is zero on the orthogonal complement to $\e_7\subset \odot^2\Co^{56}$. Thus $R\in\odot^2 \e_7$. As
the  $\e_7$-module $\odot^2 \e_7$ can be decomposed as $\odot^2 \e_7=V^{\e_7}_{2\pi_1}\oplus V^{\e_7}_{\pi_6}\oplus\Co$.
We have $\dim V^{\e_7}_{2\pi_1}=7371$ and $\dim V^{\e_7}_{\pi_6}=1539$. Not that  each such component consists of
curvature tensors if and only if it contains at least one non-zero curvature tensor (such idea was used by
D.V.~Alekseevsky in \cite{Al} to find the spaces of curvature tensors for irreducible holonomy algebras of Riemannian
manifolds).
 As the $\sl(8,\Co)$-module
$\odot^2 \e_7=\odot^2 (\sl(8,\Co)\oplus\Lambda^4\Co^8)$ can be decomposed as $\odot^2 \e_7=\odot^2
\sl(8,\Co)\oplus\odot^2\Lambda^4\Co^8\oplus W$, where $W\subset\odot^2 \e_7$ consists of symmetric maps interchanging
$\sl(8,\Co)$ and $\Lambda^4\Co^8$. An important fact is that in the decomposition of $\odot^2 \e_7$  into the direct sum
of irreducible $\sl(8,\Co)$-modules there is only one  summand of dimension greater than one that appears  twice: $\odot^2
(\sl(8,\Co))$ and $\odot^2\Lambda^4$ contain a submodule isomorphic to $V^{\sl(8,\Co)}_{\pi_2+\pi_6}$. As a consequence,
any $R\in \RR(\e_7)$ is a sum of elements of $\RR(\e_7)$ contained in $(\odot^2 \sl(8,\Co))/V^{\sl(8,\Co)}_{\pi_2+\pi_6}$,
$(\odot^2\Lambda^4\Co^8)/V^{\sl(8,\Co)}_{\pi_2+\pi_6}$, $W$ and in $2V^{\sl(8,\Co)}_{\pi_2+\pi_6}\subset\odot^2
\sl(8,\Co)\oplus\odot^2\Lambda^4\Co^8$. We have $\odot^2\Co^{56}= \odot^2(\Lambda^2\Co^8)\oplus (\Lambda^2\Co^8\otimes
\Lambda^2(\Co^8)^*)\oplus\odot^2(\Lambda^2\Co^8)=(\Lambda^2\Co^8\oplus\Lambda^4\Co^8)\oplus(V^{\sl(8,\Co)}_{\pi_2+\pi_6}\oplus\sl(8,\Co)\oplus\Co)\oplus
(\Lambda^2(\Co^8)^*\oplus\Lambda^4\Co^8),$ as the $\sl(8,\Co)$-module. This decomposition makes clear the inclusion
$\e_7\subset\odot^2\Co^{56}$ and the behavior of the elements of $W$. Suppose that $R\in\RR(\e_7)\cap W$. Let
$x,y\in\odot^2(\Lambda^2\Co^8)$. By the definition of $W$, we have $R(x,y)\in\sl(8,\Co)$. Applying the Bianchi identity to
$x,y,z\in \odot^2(\Lambda^2\Co^8)$, we get $R|_{\odot^2(\Lambda^2\Co^8)\otimes\Lambda^2\Co^8}\in \RR(\sl(8,\Co) \text{
acting on } \Lambda^2\Co^8)$. Moreover, since $R$ is symmetric, this element is zero if and only if $R=0$. From the above
we know that
 $\RR(\sl(8,\Co) \text{ acting on }\Lambda^2\Co^8)\simeq
\Lambda^2(\Co^8)^*\otimes\odot^2(\Co^8)^*\simeq V^{\sl(8,\Co)}_{\pi_5+\pi_7}\oplus V^{\sl(8,\Co)}_{\pi_6+2\pi_7}$. Thus
$\RR(\e_7)\cap W$ is isomorphic to a submodule of the $\sl(8,\Co)$-module  $V^{\sl(8,\Co)}_{\pi_5+\pi_7}\oplus
V^{\sl(8,\Co)}_{\pi_6+2\pi_7}$. Note that $\dim V^{\sl(8,\Co)}_{\pi_5+\pi_7}=378$ and $\dim
V^{\sl(8,\Co)}_{\pi_6+2\pi_7}=630$. Next, $\odot^2 (\sl(8,\Co))\simeq V^{\sl(8,\Co)}_{2\pi_1+2\pi_7}\oplus
V^{\sl(8,\Co)}_{\pi_2+\pi_6}\oplus\sl(8,\Co)\oplus\Co,$ $\odot^2\Lambda^4\Co^8\simeq V^{\sl(8,\Co)}_{\pi_4}\oplus
V^{\sl(8,\Co)}_{\pi_2+\pi_6}\oplus\Co$, and $\dim V^{\sl(8,\Co)}_{2\pi_1+2\pi_7}= 1232$,
$V^{\sl(8,\Co)}_{\pi_2+\pi_6}=720$, $\dim V^{\sl(8,\Co)}_{2\pi_4}=1764$. Analyzing the dimensions, we conclude that
$\RR(\e_7)=0$ (in other words, $W\backslash\R(\e_7)$ containes non-trivial elements from the submodules
$V^{\e_7}_{2\pi_1}$ and $V^{\e_7}_{\pi_6}\subset\odot^2 \e_7$, and as a consequence these submodules do not contain
non-zero curvature tensors). Since $\e_7\subset\sp(56,\Co)$, we get $\RR(\e_7\oplus \Co)=0$.

\vskip0.2cm

$(xix)$:  $\f_4\subset\gl(26,\Co)$ is not a skew-Berger algebra. In \cite{Adams} there is the following description of
this representation. The Lie algebra $\f_4$ admits the structure of $\mathbb{Z}_2$-graded Lie algebra:
$\f_4=\so(9,\Co)\oplus\Delta$, where $\Delta$ is the representation space for the spin representation of $\so(9,\Co)$. The
representation space $\Co^{26}$ is decomposed into the direct sum $\Co^{26}=\Co\oplus\Co^9\oplus\Delta$. The elements of
the subalgebra $\so(9,\Co)\subset\f_4$ preserve these components, annihilate $\Co$ and act naturally on $\Co^9$ and
$\Delta$. Elements of $\Delta\subset\f_4$ take $\Co$ and $\Co^9$ to $\Delta$ (multiplication by constants and the Clifford
multiplication, respectively), and take $\Delta$ to $\Co\oplus\Co^9$ (the charge conjugation plus the natural map
assigning a vector to a pair of spinors). Let $R\in\RR(\f_4)$. Decompose it as the sum $R=S+T$, where $S$ and $T$ take
values in $\so(9,\Co)$ and $\Delta$, respectively. Let $\lambda,\mu\in\Co$, $x,y,z\in\so(9,\Co)$, and $X,Y,Z\in\Delta$.
Applying the Bianchi identity to $X$, $Y$ and $Z$, we get that $S|_{\odot^2\Delta\otimes\Delta}\in\RR(\spin(9,\Co))$. On
the other hand, $\RR(\spin(9,\Co))=0$, consequently, $S(X,Y)=0$. Applying the Bianchi identity to $\lambda$, $x$ and $X$,
we get $T(\lambda,x)=0$ and $S(X,\lambda)=0$. Applying the Bianchi identity to $\lambda$, $X$ and $Y$, we get $T(X,Y)=0$.
Applying the Bianchi identity to $X$, $Y$ and $x$, we get $S(X,x)Y+S(Y,x)X=0$. This means that for each fixed $x$ the map
$S(\cdot,x):\Delta\to\so(9,\Co)$ lies in the first skew-prolongation for the representation $\spin(9,\Co)$, which is
trivial, as we already know. Consequently, $S(x,X)=0$. Writing down the Bianchi identity for other vectors, we conclude
that $R=0$. Thus, $\RR(\f_4)=0$. Since this representation is orthogonal, we have $\RR(\f_4\oplus\Co)=0$.

\vskip0.5cm

Suppose now that the semisimple part $\g_s$ is not simple, then it can be written as a direct sum $\g_s=\g_1\oplus\g_2$,
where $\g_1\subset\sl(n_1,\Co)$, $\g_2\subset\sl(n_2,\Co)$, and $V=\Co^{n_1}\otimes\Co^{n_2}$.

First suppose that  $n_1,n_2\geq 3$. By the same arguments as in \cite{Sch}, each $\g_i$ must be either $\sl(n_i,\Co)$, or
$\so(n_i,\Co)$, or $\sp(n_i,\Co)$. Suppose that $\g_1=\so(n_1,\Co)$.
Consider the subalgebra $\so(n_1,\Co)\oplus
\sl(n_2,\Co)\oplus\Co\subset\sl(n_1,\Co)\oplus\sl(n_2,\Co)\oplus\Co$. Let $R_\tau\in\RR(\sl(n_1,\Co)\oplus\sl(n_2,\Co)\oplus\Co)$
be as in Section \ref{Sskew-prolongs}, where $\tau\in V^*\otimes V^*$, $V=\Co^{n_1}\otimes\Co^{n_2}$.
Suppose that
$R_\tau\in\RR(\so(n_1,\Co)\oplus \sl(n_2,\Co)\oplus\Co)$. Then we get  $A(x_1\otimes x_2,u_1\otimes u_2)-\tr(A(x_1\otimes x_2,u_1\otimes u_2))\id_{\Co^{n_1}}\in\so(n_1,\Co)$.  Hence for all $v_1,z_1\in \Co^{n_1}$ it holds
\begin{multline*} 0=\tau(x_1,x_2,v_1,u_2)g(u_1,z_1)+\tau(u_1,u_2,v_1,x_2)g(x_1,z_1)+\tau(x_1,x_2,z_1,u_2)g(v_1,u_1)\\
+\tau(u_1,u_2,z_1,x_2)g(v_1,x_1)-\tau(x_1,x_2,u_1,u_2)g(v_1,z_1)-\tau(u_1,u_2,x_1,x_2)g(v_1,z_1).
\end{multline*}
Taking $u_1=z_1$,  $x_1$ and $v_1$ mutually orthogonal, we get that $\tau(x_1,x_2,v_1,u_2)=0$ whenever $g(x_1,v_1)=0$.
Taking $v_1=x_1$ orthogonal to $u_1=z_1$ such that $g(x_1,x_1)=g(u_1,u_1)=1$, we get
$\tau(x_1,x_2,x_1,u_2)=-\tau(z_1,u_2,z_1,x_2)$ whenever $g(x_1,z_1)=0$. In particular, $\tau(z_1,u_2,z_1,x_2)$ does not
depend on $z_1$ under the condition $g(x_1,z_1)=0$ and $g(z_1,z_1)=1$. Considering an orthonormal basis of $\Co^{n_1}$, we
conclude that $\tau(x_1,x_2,u_1,u_2)=g(x_1,u_1)w(x_2,u_2)$ for some skew-symmetric (not necessary non-degenerate) bilinear
form on $\Co^{n_2}$. Thus, $\RR(\so(n_1,\Co)\oplus \sl(n_2,\Co)\oplus\Co)=\RR(\so(n_1,\Co)\oplus \sl(n_2,\Co))\simeq
\Lambda^2(\Co^{n_2})^*$. Using this, it is easy to get that $\RR^\nabla (\so(n_1,\Co)\oplus \sl(n_2,\Co))=0$. Thus
$\so(n_1,\Co)\oplus \sl(n_2,\Co)$ is a symmetric Berger algebra, but it does not admit a curvature tensor that is
annihilated by this algebra and such that its image coincides with this Lie algebra (this representation does not appear
in Table \ref{tabsym}). Consequently, any absolutely irreducible real representation $\h\subset\sl(n_1n_2,\Real)$ with
such complexification can not appear as the holonomy algebra of a linear torsion-free connection on a purely odd
supermanifold (if $\h\subset\sl(n_1n_2,\Real)$ appears as the holonomy algebra of an odd supermanifold, then the curvature
tensor of this manifold is parallel; hence its value at the point  must be annihilated by $\h$ and the image of this
curvature tensor must coincide with $\h$).

Similarly, if $n_1$ is even, then $\RR(\sp(n_1,\Co)\oplus \sl(n_2,\Co)\oplus\Co)=\RR(\sp(n_1,\Co)\oplus
\sl(n_2,\Co))\simeq \odot^2(\Co^{n_1})^*$.  We also get that $\RR(\so(n,\Co)\oplus \sp(2m,\Co))$ is one-dimensional and it
is spanned by the curvature tensor $R_\tau$ with $\tau=g\otimes w$.

\vskip0.2cm

By the same arguments as in \cite{Sch} it can be proved that if $n_1=2$ and $\g$ is a non-symmetric Berger algebra, then
$\g_1=\sl(2,\Co)$ and $\g_2$ is one of $\gl(n_2,\Co)$, $\sl(n_2,\Co)$, $\so(n_2,\Co)$, or $\sp(n_2,\Co)$.

The theorem is proved. $\Box$

\section{An outlook to the general case}\label{Outlook}

Consider the identity representation of the orthosymplectic Lie superalgebra $\osp(n|2m,\Co)\subset\gl(n|2m,\Co)$ on the
vector superspace $\Co^n\oplus\Pi(\Co^{2m})$. Recall that $\osp(n|2m,\Co)$ is the supersubalgebra of $\gl(n|2m,\Co)$
preserving the form $g+\Omega$, where $g$ is the standard non-degenerate symmetric bilinear form on $\Co^n$ and $\Omega$
is the standard non-degenerate skew-symmetric bilinear form on $\Co^{2m}$. For the even part of $\osp(n|2m,\Co)$ we have
$\osp(n|2m,\Co)_{\bar 0}=\so(n,\Co)\oplus\sp(2m,\Co)$. Note that $\so(n,\Co)$ annihilates  $\Pi(\Co^{2m})$ and acts on
$\Co^n$ in the natural way. Similarly, $\sp(2m,\Co)$ annihilates  $\Co^{n}$ and acts on $\Pi(\Co^{2m})$ in the natural
way. It is easy to describe the space $\R(\osp(n|2m,\Co))$ using the method of \cite{GalLorCurv}. In particular, the
following curvature tensors take values in $\osp(n|2m,\Co)_{\bar 0}$: $R_0|_{\Co^n\wedge\Co^n}$,
$R_0|_{\Pi(\Co^{2m})\wedge\Pi(\Co^{2m})}$, and $R_1|_{\Co^n\otimes\Pi(\Co^{2m})}$, where $R_0\in\R(\osp(n|2m,\Co))_{\bar
0}$ and $R_1\in\R(\osp(n|2m,\Co))_{\bar 1}$.

For a subalgebra $\h\subset\so(n,\Co)$ define the space of weak curvature tensors of type $\h$,
$$\P_g(\h)=\{P\in(\Co^n)^*\otimes \h|\,g(P(x)y,z)+g(P(y)z,x)+g(P(z)x,y)=0\text{ for all } x,y,z\in\Co^n\}.$$
This space was introduced in \cite{GalLorCurv,Leistner} and it appears if one considers the space of curvature tensors for the holonomy algebras of Lorentzian manifolds. Note that if $R\in\R(\h)$, then for any fixed $x\in\Co^n$ it holds $R(\cdot,x)\in\P_g(\h)$. A subalgebra $\h\subset\so(n,\Co)$ is called {\it a weak-Berger algebra} if it is spanned by the images of the elements of $\P_g(\h)$.

The following important theorem is proved by T.~Leistner in \cite{Leistner}.

\begin{theorem}\label{Leistner} Let $\h\subset\so(n,\Co)$ be an irreducible subalgebra. Then $\h$ is a weak-Berger algebra if and only if it is a Berger algebra.\end{theorem}

Similarly, for a subalgebra $\h\subset\sp(2m,\Co)$ define the space of weak skew-curvature tensors of type $\h$,
$$P_\Omega(\h)=\{P\in(\Co^n)^*\otimes \h|\,\Omega(P(x)y,z)+\Omega(P(y)z,x)+\Omega(P(z)x,y)=0\text{ for all } x,y,z\in\Co^{2m}\}.$$
 Note that if $R\in\RR(\h)$, then for any fixed $x\in\Co^{2m}$ it holds $R(\cdot,x)\in\P_\Omega(\h)$. We call a subalgebra $\h\subset\sp(2m,\Co)$
{\it a weak-skew-Berger algebra} if it is spanned by the images of the elements of $\P_\Omega(\h)$.

In view of Theorem \ref{Leistner}, we get the following.

\begin{hyp}\label{hyp1} Let $\h\subset\sp(2m,\Co)$ be an irreducible subalgebra.
Then $\h$ is a weak-skew-Berger algebra if and only if it is a skew-Berger algebra.\end{hyp}

Let now $R_0$ and $R_1$ be as above, then $\pr_{\so(n,\Co)}\circ R_0|_{\Co^n\wedge\Co^n}\in \R(\so(n,\Co))$,
$\pr_{\sp(2m,\Co)}\circ R_0|_{\Pi(\Co^{2m})\wedge\Pi(\Co^{2n})}\in \R(\sp(2m,\Co) \text{ acting on } \Pi(\Co^{2m}))$, and
for any fixed $x_0\in \Co^n$ and  $x_1\in \Pi(\Co^{2m})$ it holds\\ $\pr_{\so(n,\Co)}\circ
R_1(\cdot,x_1)|_{\Co^n}\in\P_g(\so(n,\Co))$ and $\pr_{\sp(2m,\Co)}\circ
R_1(\cdot,x_0)|_{\Pi(\Co^{2m})}\in\P_\Omega(\sp(2m,\Co))$. On the other hand, there is no such obvious restrictions on
$\pr_{\sp(2m,\Co)}\circ R_0|_{\Co^n\wedge\Co^n}$ and $\pr_{\so(n,\Co)}\circ R_0|_{\Pi(\Co^{2m})\wedge\Pi(\Co^{2n})}$.

Suppose now that we have a simple supersubalgebra $\g\subset\osp(n|2m,\Co)$ such that the  representations of $\g_{\bar
0}$ in the both $\Co^n$ and $\Pi(\Co^{2m})$ are faithful. This is not the case for the identity representation of
$\osp(n|2m,\Co)$, but this is the case for the adjoint representations of the simple Lie superalgebras of classical type
and it seems to be the case for the most of orthosymplectic representations of simple Lie superalgebras. In this case
$\pr_{\sp(2m,\Co)}\circ R_0|_{\Co^n\wedge\Co^n}$ and $\pr_{\so(n,\Co)}\circ R_0|_{\Pi(\Co^{2m})\wedge\Pi(\Co^{2n})}$ are
determined, respectively, by $\pr_{\so(n,\Co)}\circ R_0|_{\Co^n\wedge\Co^n}$ and \\ $\pr_{\sp(2m,\Co)}\circ
R_0|_{\Pi(\Co^{2m})\wedge\Pi(\Co^{2n})}$,
 and therefore are strongly restricted.

Module Hypothesis \ref{hyp1} we get.

\begin{theorem} Let $\g\subset\osp(n|2m,\Co)$ be a simple subalgebra such that the  representations of $\g_{\bar 0}$ in  both $\Co^n$ and $\Pi(\Co^{2m})$ are fatheful. Then $\g_{\bar 0}\subset \so(n,\Co)$ is a Berger algebra and $\g_{\bar 0}\subset \sp(2n,\Co)$ is a skew-Berger algebra.
\end{theorem}

Similarly, for a simple Berger subalgebra $\g\subset\gl(n|m,\Co)$  such that the  representations of $\g_{\bar 0}$ in
both $\Co^n$ and $\Pi(\Co^{m})$ are faithful, we expect that there are two ideals $\g_1,\g_2\in\g_{\bar 0}$ such that
$\g_1+\g_2=\g_{\bar 0}$ and $\g_1\subset \gl(n,\Co)$ is a Berger algebra, and $\g_{2}\subset \gl(m,\Co)$ is a skew-Berger
algebra.

In another paper we will discuss the ideas of this section in details.

\bibliographystyle{unsrt}

\vskip1cm

Department of Algebra and Geometry, Masaryk University in Brno, Kotl\'a\v rsk\'a~2, 611~37 Brno, Czech Republic

{\it E-mail address}: galaev@math.muni.cz

\end{document}